\documentclass[11pt]{article}

\usepackage{xcolor}
\definecolor{labelkey}{rgb}{0,0.08,0.45}
\definecolor{refkey}{rgb}{0,0.6,0.0}
\definecolor{Brown}{rgb}{0.45,0.0,0.05}
\definecolor{dgreen}{rgb}{0.00,0.40,0.00}
\definecolor{dblue}{rgb}{0,0.08,0.45}

\usepackage{amsmath}
\usepackage{amssymb}
\usepackage{theorem}
\newcommand{\emp}{\ensuremath{\varnothing}}
\newcommand{\infconv}{\ensuremath{\mbox{\small$\,\square\,$}}}

\newcommand{\envv}[2]{\mathrm{env}_{#1}({#2})}

\newcommand{\scal}[2]{\left.\left\langle{#1}\:\right| {#2}  \right\rangle} 
\newcommand{\abscal}[2]{\left|\left\langle{{#1}\mid{#2}}%
\right\rangle\right|} 
\newcommand{\HH}{\ensuremath{\mathcal H}}

\newcommand{\RR}{\ensuremath{\mathbb R}}

\newcommand{\RPP}{\ensuremath{\,\left]0,+\infty\right[}}
\newcommand{\RX}{\ensuremath{\,\left]-\infty,+\infty\right]}}
\newcommand{\NN}{\ensuremath{\mathbb N}}
\newcommand{\IDD}{\ensuremath{\operatorname{int}\operatorname{dom}\varphi}}

\newcommand{\dom}{\ensuremath{\operatorname{dom}}}

\newcommand{\prox}{\ensuremath{\operatorname{Prox}}}
\newcommand{\intdom}{\ensuremath{\operatorname{int}\operatorname{dom}}\,}

\newcommand{\Nf}{\ensuremath{\nabla f}}

\newcommand{\Id}{\ensuremath{\operatorname{Id}}}

\newcommand{\pinf}{\ensuremath{+\infty}}
\newcommand{\menge}[2]{\big\{{#1} \mid {#2}\big\}}

\newtheorem{theorem}{Theorem}[section]

\newtheorem{corollary}[theorem]{Corollary}
\newtheorem{proposition}[theorem]{Proposition}

\theoremstyle{plain}{\theorembodyfont{\rmfamily}
}
\theoremstyle{plain}{\theorembodyfont{\rmfamily}
}
\theoremstyle{plain}{\theorembodyfont{\rmfamily}
}
\theoremstyle{plain}{\theorembodyfont{\rmfamily}
}

\theoremstyle{plain}{\theorembodyfont{\rmfamily}
\newtheorem{remark}[theorem]{Remark}}

\begin{document}

\title{\bfseries The Baillon-Haddad Theorem Revisited}

\author{{\bfseries Heinz H. Bauschke}\thanks{Research partially 
supported by the Natural Sciences and Engineering Research Council of 
Canada and by the Canada Research Chair Program.}\\[1mm]
{\em Mathematics, University of British Columbia Okanagan}\\
{\em Kelowna, B.C.~V1V 1V7, Canada}\\
{\em heinz.bauschke@ubc.ca}\\[2mm]
\and
{\bfseries Patrick L. Combettes}\thanks{Research supported by the 
Agence Nationale de la Recherche under grant
ANR-08-BLAN-0294-02.}\\[1mm]
{\em UPMC Universit\'e Paris 06}\\
{\em Laboratoire Jacques-Louis Lions -- UMR CNRS 7598}\\
{\em 75005 Paris, France}\\
{\em plc@math.jussieu.fr}}

\date{~}

\maketitle

\vskip 8mm

\begin{abstract}
\noindent
In 1977, Baillon and Haddad proved that if the gradient of a convex 
and continuously differentiable function is nonexpansive, then it 
is actually firmly nonexpansive. This result, which has become known 
as the Baillon-Haddad theorem, has found many applications in 
optimization and numerical functional analysis. In this note, we
propose short alternative proofs of this result and strengthen its
conclusion.

\end{abstract}

\noindent
\textbf{2000 Mathematics Subject Classification:}\\
Primary 47H09, 90C25; Secondary 26A51, 26B25, 46C05, 47H05, 52A41. 

\noindent
\textbf{Key Words:}\\
Backward-backward splitting,
Bregman distance,
cocoercivity,
convex function, 
Dunn property, 
firmly nonexpansive,
forward-backward splitting,
gradient, 
inverse strongly monotone,
Moreau envelope,
nonexpansive,
proximal mapping,
proximity operator.

\section{Introduction}
Throughout, $\HH$ is a real Hilbert space with scalar product 
$\scal{\cdot}{\cdot}$ and induced norm $\|\cdot\|$. 
Let $C$ be a nonempty subset of $\HH$, let $T\colon C\to\HH$, and let 
$\beta\in\RPP$. Then $T$ is $1/\beta$-cocoercive if
(this property is also known as the Dunn property
or inverse strong monotonicity) 
\begin{equation}
\label{e:rio-mai2009-1}
(\forall x\in C)(\forall y\in C)\quad
\beta\scal{x-y}{Tx-Tx}\geq\|Tx-Ty\|^2,
\end{equation}
and $T$ is $\beta$-Lipschitz continuous if
\begin{equation}
\label{e:rio-mai2009-2}
(\forall x\in C)(\forall y\in C)\quad
\|Tx-Tx\|\leq\beta\|x-y\|^2.
\end{equation}
When $\beta=1$, \eqref{e:rio-mai2009-1} means that $T$ is 
firmly nonexpansive and \eqref{e:rio-mai2009-2} that $T$ is nonexpansive.
Cocoercivity arises in various areas of 
optimization and nonlinear analysis, e.g., 
\cite{Atto10,Byrn08,Opti04,Dunn76,Lion79,Liuf98,Tsen91,Zhud96}. It 
follows from the Cauchy-Schwarz inequality that 
$1/\beta$-cocoercivity
implies $\beta$-Lipschitz continuity. However, the converse fails; take 
for instance $T=-\Id$, which is nonexpansive but not firmly
nonexpansive. In 1977, Baillon and Haddad showed that, if $C=\HH$ and
$T$ is the gradient of a convex function, then 
\eqref{e:rio-mai2009-1} and \eqref{e:rio-mai2009-2} coincide.
This remarkable result, which has important applications in 
optimization (see for instance \cite{Smms05,Yama04}), has become known 
as the Baillon-Haddad theorem.

\begin{theorem}[Baillon-Haddad]
\emph{\cite[Corollaire~10]{Bail77}}
\label{t:B-H}
Let $f\colon\HH\to\RR$ be convex, Fr\'echet differentiable on $\HH$, 
and such that $\nabla f$ is $\beta$-Lipschitz continuous for some
$\beta\in\RPP$. Then $\nabla f$ is $1/\beta$-cocoercive. 
\end{theorem}

In \cite{Bail77}, Theorem~\ref{t:B-H} was derived from a more general
result concerning $n$-cyclically monotone operators in normed vector
spaces. Since then, direct proofs have been proposed, such as 
\cite[Lemma~6.7]{GolTre}, \cite[Theorem~X.4.2.2]{HULL2},
and \cite[Proposition~12.60]{Rock98} for Euclidean spaces.
These approaches rely on convex analytical and integration 
arguments. An infinite dimensional proof can be found in 
\cite[Remark~3.5.2]{Zali02}, as a corollary to results on 
the properties of uniformly smooth convex functions. 

The goal of our paper is to provide new insights into the
Baillon-Haddad theorem. In Section~\ref{sec:2}, we propose a short 
new proof of Theorem~\ref{t:B-H} and present additional equivalent
conditions, thus making a connection with lesser known parts of
Moreau's classical paper \cite{More65}. In Section~\ref{sec:3}, we
provide a second order variant of the Baillon-Haddad theorem that
partially extends work by Dunn \cite{Dunn76}.

{\bfseries Notation and background.} 
Our notation is standard: $\Gamma_0(\HH)$ is the class of proper lower
semicontinuous convex functions from $\HH$ to $\RX$ and $\infconv$ 
denotes infimal convolution. The conjugate of a function
$f\colon\HH\to\RX$ is denoted by $f^*$, and its subdifferential 
by $\partial f$.
For background on convex analysis, we refer the reader to
\cite{HULL2,Rock70,Zali02}.

\section{An enhanced Baillon-Haddad theorem}
\label{sec:2}

Let us start with some standard facts on Moreau envelopes and 
proximity operators; we refer the reader to Moreau's original 
paper \cite{More65} and to \cite{Atto84,Smms05,Rock98} for details and 
complements. Let $\varphi\in\Gamma_0(\HH)$ and let $\gamma\in\RPP$. 
The Moreau envelope of $\varphi$ of index $\gamma$ is the finite 
continuous convex function 
\begin{equation}
\label{e:moreau1}
\envv{\gamma}{\varphi}=\varphi\infconv\big(q/\gamma\big),
\quad\text{where}\quad
q=\text{\footnotesize$\frac12$}\|\cdot\|^2.
\end{equation}
Moreau's decomposition asserts that 
\begin{equation}
\label{e:moreau2}
\envv{1/\gamma}{\varphi}+\envv{\gamma}{\varphi^*}\circ\,(\gamma\Id)
=\gamma q.
\end{equation}
The proximity operator (or proximal mapping) of $f$ is the operator 
$\prox_\varphi=(\Id+\partial\varphi)^{-1}$; it
maps each $x\in\HH$ to the unique minimizer of 
the function $y\mapsto\varphi(y)+q(x-y)$. The Moreau envelope
$\envv{1}{\varphi}$ is Fr\'echet differentiable with gradient 
$\nabla\envv{1}{\varphi}=\prox_{\varphi^*}$. Hence, \eqref{e:moreau2}
yields
\begin{equation}
\label{e:moreau4}
\nabla\,\envv{1/\gamma}{\varphi}
=\prox_{\gamma\varphi^*}\circ\,(\gamma\Id)
=\gamma(\Id-\prox_{\varphi/\gamma}).
\end{equation}
Moreover,
\begin{equation}
\label{e:moreau3}
\prox_\varphi\colon\HH\to\HH\;\;\text{is firmly nonexpansive.}
\end{equation}

We are now ready to present the main result of this section, which
strengthens the conclusion of Theorem~\ref{t:B-H} by providing four
additional equivalent conditions and a short new proof.

\begin{theorem}
\label{t:rio2009-05-26}
Let $f\in\Gamma_0(\HH)$, let $\beta\in\RPP$, and set $h=f^*-q/\beta$. 
Then the following are equivalent.
\begin{enumerate}
\item
\label{t:rio2009-05-26i}
$f$ is Fr\'echet differentiable on $\HH$ and 
$\Nf$ is $\beta$-Lipschitz continuous.
\item
\label{t:rio2009-05-26ii}
$\beta q-f$ is convex.
\item
\label{t:rio2009-05-26iii}
$f^*-q/\beta$ is convex (i.e., $f^*$ is $1/\beta$-strongly convex).
\item
\label{t:rio2009-05-26iv}
$h\in\Gamma_0(\HH)$ and
$f=\envv{1/\beta}{h^*}=\beta q-\envv{\beta}{h}\circ\,\beta\Id$.
\item
\label{t:rio2009-05-26v}
$h\in\Gamma_0(\HH)$ and
$\nabla f=\prox_{\beta h}\circ\,\beta\Id=\beta(\Id-\prox_{h^*/\beta})$.
\item
\label{t:rio2009-05-26vi}
$f$ is Fr\'echet differentiable on $\HH$ and 
$\nabla f$ is $1/\beta$-cocoercive.
\end{enumerate}
\end{theorem}
\begin{proof}
\ref{t:rio2009-05-26i}$\Rightarrow$\ref{t:rio2009-05-26ii}:
By Cauchy-Schwarz, $(\forall x\in\HH)(\forall y\in\HH)$ 
$\scal{x-y}{\beta x-\nabla f(x)-\beta y+
\nabla f(y)}=\beta\|x-y\|^2-\scal{x-y}{\nabla
f(x)-\nabla f(y)}\geq\|x-y\|(\beta\|x-y\|-\|\nabla f(x)-\nabla f(y)\|)
\geq 0$. Hence, $\nabla(\beta q-f)=\beta\Id-\nabla f$ is monotone and 
it follows that $\beta q-f$ is convex (see, e.g., 
\cite[Theorem~2.1.11]{Zali02}).

\ref{t:rio2009-05-26ii}$\Rightarrow$\ref{t:rio2009-05-26iii}:
Set $g=\beta q-f$. Then $g\in\Gamma_0(\HH)$ and therefore
$g=g^{**}$. Accordingly, 
\begin{equation}
f=\beta q-g=\beta q-g^{**}
=\beta q-\sup_{u\in\HH}\big(\scal{\cdot}{u}-g^*(u)\big)
=\inf_{u\in\HH}\big(\beta q-\scal{\cdot}{u}+g^*(u)\big).
\end{equation}
Hence
\begin{multline}
f^*=\sup_{u\in\HH}\big(\beta q-\scal{\cdot}{u}+g^*(u)\big)^*=
\sup_{u\in\HH}\big(\big(\beta q-\scal{\cdot}{u}\big)^*-g^*(u)\big)\\
=\sup_{u\in\HH}\big(q(\cdot+u)/\beta-g^*(u)\big)
=q/\beta+\sup_{u\in\HH}\big(
(\scal{\cdot}{u}+q(u))/\beta-g^*(u)\big),
\end{multline}
where the last term is convex as a supremum of affine functions.
Thus, $h$ is convex.

\ref{t:rio2009-05-26iii}$\Rightarrow$\ref{t:rio2009-05-26iv}:
Since $f\in\Gamma_0(\HH)$ and $h$ is convex, we have
$h\in\Gamma_0(\HH)$, $h^*\in\Gamma_0(\HH)$, and 
$f=f^{**}=(h+q/\beta)^*=h^*\infconv\beta q
=\envv{1/\beta}{h^*}=\beta q-\envv{\beta}{h}\circ\,\beta\Id$, 
where the last identity follows from \eqref{e:moreau2}. 

\ref{t:rio2009-05-26iv}$\Rightarrow$\ref{t:rio2009-05-26v}:
Use \eqref{e:moreau4}.

\ref{t:rio2009-05-26v}$\Rightarrow$\ref{t:rio2009-05-26vi}:
By \eqref{e:moreau3}, $\prox_{\beta h}$ is firmly nonexpansive.
Hence, it follows from \eqref{e:rio-mai2009-1} that 
$\nabla f=\prox_{\beta h}\circ\,\beta\Id$ 
is $1/\beta$-cocoercive.

\ref{t:rio2009-05-26vi}$\Rightarrow$\ref{t:rio2009-05-26i}:
Apply the Cauchy-Schwarz inequality.
\end{proof}

\begin{remark}
Some comments regarding Theorem~\ref{t:rio2009-05-26} are in order.
\begin{itemize}
\item[(a)]
The proof of the implication
\ref{t:rio2009-05-26i}$\Rightarrow$\ref{t:rio2009-05-26vi}, i.e., of 
the Baillon-Haddad theorem (Theorem~\ref{t:B-H}) appears to be new and 
shorter than those found in the literature. In addition,
Theorem~\ref{t:rio2009-05-26} brings to light various characterizations
of the Lipschitz continuity of the gradient of a convex function.
The equivalences 
\ref{t:rio2009-05-26ii}$\Leftrightarrow$\ref{t:rio2009-05-26iii}%
$\Leftrightarrow$\ref{t:rio2009-05-26iv} are due to Moreau, who
established them (for $\beta=1$) in \cite[Proposition~9.b]{More65} 
(see also \cite[Corollary~3]{Jbhu89}).
On the other hand, the equivalences 
\ref{t:rio2009-05-26i}$\Leftrightarrow$\ref{t:rio2009-05-26iii}%
$\Leftrightarrow$\ref{t:rio2009-05-26iv}$\Leftrightarrow$%
\ref{t:rio2009-05-26vi} are shown in Euclidean spaces in
\cite[Proposition~12.60]{Rock98} with different techniques.
\item[(b)]
Set $\beta=1$. The conclusion of Theorem~\ref{t:B-H} is that 
$\nabla f\colon\HH\to\HH$ is firmly nonexpansive. Hence, since the class 
of firmly nonexpansive operators with domain $\HH$ coincides with that 
of resolvents of maximal monotone operators \cite[Section~1.11]{Goeb84}, 
we have $\nabla f=(\Id+A)^{-1}$, for some maximal monotone operator 
$A\colon\HH\to 2^{\HH}$. However, \ref{t:rio2009-05-26v} more precisely
reveals $\nabla f$ to be the proximity operator of $h$, i.e., 
$A=\partial h=\partial f^*-\Id$.
\item[(c)]
Let $f_1\in\Gamma_0(\HH)$, let $f_2\colon\HH\to\RR$ be convex and 
differentiable with a Lipschitz continuous gradient, 
and consider the problem of minimizing 
$f_1+f_2$. Without loss of generality (rescale), we assume that
the Lipschitz constant of $\nabla f_2$ is $\beta=1$. A standard algorithm 
for solving this problem is the forward-backward algorithm 
\cite{Smms05,Tsen91}
\begin{equation}
\label{e:main1}
x_0\in\HH\quad\text{and}\quad(\forall n\in\NN)\quad
x_{n+1}=\prox_{\gamma_n f_1}\big(x_n-\gamma_n\nabla f_2(x_n)\big),
\quad 0<\gamma_n<2.
\end{equation}
Now set $h_2=f_2^*-q$. Then it follows from the implication
\ref{t:rio2009-05-26i}$\Rightarrow$\ref{t:rio2009-05-26v}
that $\nabla f_2=\Id-\prox_{h_2^*}$. Hence, we can
rewrite \eqref{e:main1} as
\begin{equation}
\label{e:main2}
x_0\in\HH\quad\text{and}\quad(\forall n\in\NN)\quad
x_{n+1}=\prox_{\gamma_n f_1}\big((1-\gamma_n)x_n+
\gamma_n\prox_{h_2^*}x_n\big),
\quad 0<\gamma_n<2.
\end{equation}
This shows that the forward-backward algorithm \eqref{e:main1}
is actually a backward-backward algorithm. In particular, for 
$\gamma_n\equiv 1$, we recover the basic backward-backward iteration
$x_{n+1}=\prox_{f_1}\prox_{h_2^*}x_n$.
\end{itemize}
\end{remark}

We conclude this section with an alternative formulation of the
Baillon-Haddad theorem that brings into play Bregman distances.
Recall that if $\varphi\in\Gamma_0(\HH)$ is G\^ateaux differentiable on 
$\IDD\neq\varnothing$, the associated Bregman distance $D_\varphi$ is 
defined by
\begin{equation}
\label{e:sicon2003}
D_\varphi\colon\HH\times\HH\to\left[0,\pinf\right]\colon (x,y) \mapsto 
\begin{cases}
\varphi(x)-\varphi(y)-\scal{x-y}{\nabla\varphi(y)}, 
&\text{if}\;\;y\in\IDD;\\
\pinf, & \text{otherwise.}
\end{cases}
\end{equation}

\begin{corollary}
\label{p:main3}
Let $\beta\in\RPP$, and let $f\colon\HH\to\RR$ be convex, Fr\'echet 
differentiable on $\HH$, and such that $f^*$ is G\^ateaux differentiable 
on $\intdom f^*\neq\emp$. Then the following are equivalent. 
\begin{enumerate}
\item 
\label{c:main3i} 
$\Nf$ is $\beta$-Lipschitz continuous.
\item 
\label{c:main3ii}
$(\forall x\in\HH)(\forall y\in\HH)$ $D_f(x,y)\leq\beta q(x-y)$.
\item 
\label{c:main3iii}
$(\forall x^*\in \HH)(\forall y^*\in \HH)$
$\beta D_{f^*}(x^*,y^*)\geq q(x^*-y^*)$.
\end{enumerate}
\end{corollary}
\begin{proof}
\ref{c:main3i}$\Leftrightarrow$\ref{c:main3ii}:
Set $g=\beta q-f$. Then $g$ is Fr\'echet differentiable on $\dom g=\HH$ 
and $\nabla g=\beta\Id-\nabla f$. Hence, it follows from the equivalence 
\ref{t:rio2009-05-26i}$\Leftrightarrow$\ref{t:rio2009-05-26ii} in 
Theorem~\ref{t:rio2009-05-26} and \eqref{e:sicon2003} that \ref{c:main3i} 
$\Leftrightarrow$ $g\in\Gamma_0(\HH)$ is Fr\'echet differentiable on 
$\intdom f=\HH$ $\Leftrightarrow$ $(\forall x\in\HH)(\forall y\in\HH)$ 
$D_g(x,y)\geq 0$ $\Leftrightarrow$ $(\forall x\in\HH)(\forall y\in\HH)$ 
$D_f(x,y)\leq\beta q(x-y)$.

\ref{c:main3i}$\Leftrightarrow$\ref{c:main3iii}:
Set $h=f^*-q/\beta$. Then $h$ is G\^ateaux differentiable on 
$\intdom h=\intdom f^*$, with $\nabla h=\nabla f^*-(1/\beta)\Id$.
Hence, in view of the equivalence
\ref{t:rio2009-05-26i}$\Leftrightarrow$\ref{t:rio2009-05-26iii} in 
Theorem~\ref{t:rio2009-05-26} and \eqref{e:sicon2003}, 
\ref{c:main3i} $\Leftrightarrow$ $h\in\Gamma_0(\HH)$ is 
G\^ateaux differentiable on $\intdom h=\intdom f^*$ $\Leftrightarrow$ 
$(\forall x^*\in \HH)(\forall y^*\in\HH)$ 
$D_h(x^*,y^*)\geq 0$ $\Leftrightarrow$ 
$(\forall x^*\in \HH)(\forall y^*\in \HH)$ 
$D_{f^*}(x^*,y^*)\geq q(x^*-y^*)/\beta$.
\end{proof}

\section{A second order Baillon-Haddad theorem}
\label{sec:3}

Under the more restrictive assumption that 
the underlying convex function is twice continuously differentiable, 
we shall obtain in Theorem~\ref{t:main2} a very short and transparent 
proof inspired by the work of Dunn \cite{Dunn76}. 
We require two preliminary propositions.

\begin{proposition}  \label{p:charnonexp}
Let $C$ be a nonempty open convex subset of $\HH$, let ${\mathcal B}$ 
be a real Banach space, and let $G\colon C\to {\mathcal B}$ be 
continuously Fr\'echet differentiable on $C$. 
Then $G$ is nonexpansive if and only if $(\forall x\in C)$
$\|\nabla G(x)\| \leq 1$.
\end{proposition}
\begin{proof}
Let $x\in C$ and let $y\in \HH$. Suppose that $G$ is nonexpansive. 
For every $t\in\RPP$ sufficiently small, $x+ty \in C$ and
hence $\|G(x+ty)-G(x)\|/t \leq \|y\|$. 
Letting $t\downarrow 0$, we deduce that $\|(\nabla G(x))y\| \leq \|y\|$. 
Since $y$ was chosen arbitrarily, we conclude that 
$\|\nabla G(x)\| \leq 1$.
Conversely, if $y\in C$, we derive from the mean value theorem 
(see, e.g., \cite[Theorem~5.1.12]{Denk03}) that
$\|G(y)-G(x)\|\leq\|y-x\|\sup_{z\in [x,y]}\|\nabla G(z)\|\leq\|y-x\|$. 
\end{proof}

Let $A\colon\HH\to\HH$ and $B\colon\HH\to\HH$ be self-adjoint 
bounded linear operators. Then $A$ is positive, written $A \succeq 0$, if 
$(\forall x\in\HH)$ $\scal{x}{Ax}\geq 0$. We write $A\succeq B$ if 
$A-B\succeq 0$. The following result is part of the folklore.

\begin{proposition} \label{p:squeeze}
Let $A\colon\HH\to\HH$ be a bounded self-adjoint linear operator.
Then $\|A\|\leq 1$ if and only if $\Id\succeq A\succeq -\Id$.
\end{proposition}
\begin{proof}
Assume that $\HH\neq\{0\}$ and set $S=\menge{x\in \HH}{\|x\|=1}$.
Then $\Id\succeq A$ $\Leftrightarrow$
$(\forall x \in \HH)$ 
$\scal{x}{x}\geq\scal{x}{Ax}$ 
$\Leftrightarrow$
$(\forall x \in S)$ 
$1= \scal{x}{x}\geq\scal{x}{Ax}$. 
Similarly, $A \succeq -\Id$ $\Leftrightarrow$
$(\forall x \in S)$ $\scal{x}{Ax}\geq -1$.
Hence $\Id \succeq A \succeq -\Id$ $\Leftrightarrow$
$(\forall x \in S)$ $\abscal{x}{Ax} \leq 1$
$\Leftrightarrow$
$\|A\|=\sup_{x \in S}$ $\abscal{x}{Ax}\leq 1$.
\end{proof}

The main result of this section is a Baillon-Haddad theorem for 
twice continuously Fr\'echet differentiable convex functions.
It extends \cite[Theorem~4]{Dunn76}, which
assumed in addition that $f$ has full domain and uniformly bounded
Hessians. 

\begin{theorem} 
\label{t:main2}
Let $C$ be a nonempty open convex subset of $\HH$, let $f\colon C\to\RR$ 
be convex and twice continuously Fr\'echet differentiable on $C$, and let 
$\beta\in\RPP$. Then $\nabla f$ is $\beta$-Lipschitz continuous if 
and only if it is $1/\beta$-cocoercive. 
\end{theorem}
\begin{proof}
Define two operators on $C$ by $G =(1/\beta)\Nf$ and by 
$H=\nabla G =(1/\beta)\nabla^2 f$. Under our assumptions, the convexity
of $f$ is characterized by \cite[Theorem~2.1.11]{Zali02}
\begin{equation}
\label{e:2nd}
(\forall x\in\HH)\quad H(x)\succeq 0. 
\end{equation}
Hence, 
\begin{align}
\text{$\nabla f$ is $\beta$-Lipschitz continuous} &\Leftrightarrow 
\text{$G$ is nonexpansive} &&\notag \\
&\Leftrightarrow 
(\forall x \in C)\quad \|H(x)\| \leq 1  && 
\text{(by Proposition~\ref{p:charnonexp})}\notag \\
&\Leftrightarrow 
(\forall x \in C)\quad -\Id \preceq H(x) \preceq \Id &&
\text{(by Proposition~\ref{p:squeeze})}\notag\\
&\Leftrightarrow
(\forall x \in C)\quad 0\preceq H(x) \preceq \Id && 
\text{(by \eqref{e:2nd}})\notag\\
&\Leftrightarrow (\forall x \in C)\quad -\Id \preceq 2H(x)-\Id \preceq \Id
&& \notag\\
&\Leftrightarrow (\forall x \in C)\quad \|2H(x)-\Id\| \leq 1 &&
\text{(by Proposition~\ref{p:squeeze})}\notag \\
&\Leftrightarrow \text{$2G-\Id$ is nonexpansive} &&
\text{(by Proposition~\ref{p:charnonexp})}\notag \\
&\Leftrightarrow \text{$G$ is firmly nonexpansive} && 
\text{(by {\cite[Lemma~1.11.1]{Goeb84}})}\notag\\
&\Leftrightarrow \text{$\nabla f$ is $1/\beta$-cocoercive,} && \notag
\text{(by \eqref{e:rio-mai2009-1}})\notag
\end{align}
and we obtain the conclusion.
\end{proof}

In linear functional analysis, the following property is usually 
obtained via spectral theory.  
\begin{corollary}
\label{c:rio2009-05-31}
Let $A\colon\HH\to\HH$ be a positive self-adjoint bounded linear 
operator. Then $(\forall x\in\HH)$ $\|A\|\scal{x}{Ax}\geq\|Ax\|^2$. 
\end{corollary}
\begin{proof}
This is an application of Theorem~\ref{t:main2} with 
$f\colon\HH\to\RR\colon x\mapsto\scal{x}{Ax}/2$. Indeed, $f$ is 
twice continuously Fr\'echet differentiable on $\HH$ with 
$\nabla f=A$, which is $\|A\|$-Lipschitz continuous.
\end{proof}

\begin{remark}
It would be interesting to see whether 
Theorem~\ref{t:main2} holds true when the second-order assumption
is replaced by Fr\'echet differentiability. 
However, the natural approach by approximation does not appear
to be applicable; 
see \cite[Section~5]{BorNol} for pertinent comments.
\end{remark}

\end{document}